\documentclass[12pt]{article}
\usepackage[latin1]{inputenc}
\usepackage{amsmath}
\usepackage{amsfonts}
\usepackage{amssymb}
\usepackage{amsthm} 
\usepackage{euscript}

\newcommand{\eqsepv}{\; , \enspace}       
\newcommand{\eqfinv}{\; ,}                
\newcommand{\eqfinp}{\; .}

\newcommand{\mtext}[1]{\,\mbox{#1}\,} %text in maths

\newcommand{\defset}[2]{\left\{#1\:\left|\:#2\right.\right\}}
\newcommand{\sequence}[2]{\left\{#1\right\}_{#2}}           % Suite

\newcommand{\np}[1]{(#1)}                                   % Parenth\`{e}se normal
\newcommand{\bp}[1]{\big(#1\big)}                           % Parenth\`{e}se big
\newcommand{\Bp}[1]{\Big(#1\Big)}                           % Parenth\`{e}se Big

\newcommand{\bc}[1]{\big[#1\big]}                           % Crochet big
\newcommand{\RR}{{\mathbb R}} %"ensemble des reels"   
\newcommand{\NN}{{\mathbb N}} %"ensemble des entiers naturels"  
\newcommand{\EE}{{\mathbb E}} %"symbole d'esperance mathematique"
\newcommand{\FF}{{\mathbb F}} 
\newcommand{\GG}{{\mathbb G}} 
\newcommand{\TT}{{\mathbb T}} 
\newcommand{\PP}{{\mathbb P}} %"symbole de probabilite"
\newcommand{\YY}{{\mathbb Y}} 

\newcommand{\control}{u}
\newcommand{\Control}{U}
\newcommand{\CONTROL}{{\mathbb U}} 

\newcommand{\state}{x}
\newcommand{\State}{X}
\newcommand{\STATE}{{\mathbb X}}

\newcommand{\uncertain}{w}
\newcommand{\Uncertain}{W}
\newcommand{\UNCERTAIN}{{\mathbb W}} 

\newcommand{\history}{h}

\newcommand{\horizon}{T}
\newcommand{\tinitial}{t_0}

\newcommand{\dynamics}{f}
\newcommand{\Dynamics}{F}

\newcommand{\policy}{\lambda}
\newcommand{\POLICY}{\Lambda}
\newcommand{\strategy}{\policy}
\newcommand{\STRATEGY}{\POLICY}

\newcommand{\sdo}{{\mathbb A}} % state domain

\newcommand{\va}[1]{\mathbf{#1}}
\newcommand{\tribu}[1]{\EuScript{#1}}                        % Tribu
\newcommand{\1}{{\mathbf 1}}

\newcommand{\segment}[2]{#1\;\!\!:\;\!\!#2}
\newcommand{\Adapted}[2]{{\mathbb L}^0\bp{#1,#2}}
\newcommand{\Acceptable}{{\mathcal A}}
\newcommand{\ResilientStrategies}{{\STRATEGY}^{R}}
\newcommand{\ResilientStates}{{\STATE}^{R}}
\newcommand{\TimeStep}{\Delta\;\!\!t}

\title{A Mathematical Framework for Resilience:\\
Dynamics, Uncertainties, Strategies \\ and Recovery Regimes}

\def\email#1{#1}
\def\keywords#1{\textbf{Keywords:} #1}

\author{Michel De Lara\footnote{%
Universit\'e Paris-Est, Cermics (ENPC), 
F-77455 Marne-la-Vall\'ee, %France,
\email{delara@cermics.enpc.fr}}   
}

\date{\today}

\begin{document}
\maketitle

\begin{abstract}
Resilience is a rehashed concept in natural hazard management --- 
resilience of cities to earthquakes, to floods, to fire, etc. 
In a word, a system is said to be resilient if there exists a strategy that can
drive the system state back to ``normal'' after any perturbation.
What formal flesh can we put on such a malleable notion?
We propose to frame the concept of resilience in the mathematical garbs 
of control theory under uncertainty.
Our setting covers dynamical systems both in discrete or continuous time, 
deterministic or subject to uncertainties.
We will say that a system state is resilient if there exists an adaptive strategy
such that the generated state and control paths, contingent on uncertainties, 
lay within an acceptable domain of random processes, called recovery regimes.
We point out how such recovery regimes can be delineated thanks to so called
risk measures, making the connection with resilience indicators.
Our definition of resilience extends others, be they ``\`a la Holling'' 
or rooted in viability theory. Indeed, our definition of resilience 
is a form of controlability for whole random processes (regimes), 
whereas others require that the state values must belong to an 
acceptable subset of the state set.
\end{abstract}

\keywords{resilience; control theory; uncertainty; risk measures; recovery regimes}

% \pagebreak 
% \tableofcontents
% \pagebreak 

 \section{Introduction}
\label{Introduction}

Consider a system whose state evolves with time, being subject to a dynamics
driven both by controls and by external perturbations.
%Suddenly, those perturbations jump out of their normal range and induce a shock. 
The system is said to be resilient if there exists a strategy that can
drive the system state towards a normal regime, whatever the perturbations.
Basic references are
\cite{Holling:73,Martin:2004,Martin-Deffuant-Calabrese:2011,Rouge-Mathias-Deffuant:2013,Arnoldi-Loreau-Haegeman:2016}. 

In the case of fisheries, the state can be a vector of abundances at ages
of one or several species; the control can be fishing efforts; 
the external perturbations can affect mortality rates or birth functions
appearing in the dynamics (an extreme perturbation could be an El Ni\~no event,
affecting the populations renewal). 
In the case of a city exposed to earthquakes, floods or other climatic events,
the state can be a vector of capital stocks (energy reserves, energy
production units, water treatment plants, health units, etc.);
the controls would be the different investments in capital 
as well as current operations (flows in and out capital stocks);
the dynamics would express the changes in the stocks due to investment
and to day to day operations;
external perturbations (rain, wind, climatic events, etc.) would affect the
stocks by reducing them, possibly down to zero.
% shocks could be extremely severe perturbations (earthquakes, floods, etc.)
% during a limited time interval.

In Sect.~\ref{Ingredients_for_an_abstract_control_system_with_uncertainties},
we introduce basic ingredients from the
mathematical framework of control theory under uncertainty.
Thus equipped, we frame the concept of resilience in mathematical garbs 
in Sect.~\ref{Resilience:_A_Mathematical_Framework}. 
Then, in Sect.~\ref{Illustrations}, we provide illustrations of the 
abstract general framework and compare our approach with others,
``\`a la Holling'' or  the stochastic viability theory approach 
to resilience.
In Sect.~\ref{Resilience_and_risk_control/minimization},
we sketch how concepts from risk measures (introduced initially in
mathematical finance) can be imported to tackle resilience issues.
Finally, we discuss pros and cons of our approach to resilience in Sect.~\ref{sec:conclusions}.

\section{Ingredients for an abstract control system with uncertainties}
\label{Ingredients_for_an_abstract_control_system_with_uncertainties}

We outline the mathematical formulations of time, 
controls, states, Nature (uncertainties), flow (dynamics) and strategies.
As the reference \cite{Rouge-Mathias-Deffuant:2013} is the more
mathematically driven paper on resilience, we will systematically 
emphasize in what our approach differs from that of 
Roug\'e, Mathias and Deffuant.

\subsection{Time, states, controls, Nature and flow}

We lay out the basic ingredients of control theory: 
time, states, controls, Nature (uncertainties) and flow (dynamics).

\subsubsection{Time}
%\paragraph{Time.}

The \emph{time set}~$\TT$ is a (nonempty) subset of the real line~$\RR$.
The set~$\TT$ holds a minimal element~$\tinitial \in \TT$ and 
an upper bound~$\horizon$, 
which is either a maximal element when $\horizon <+\infty$
(that is, $\horizon \in \TT$)
or not when $\horizon =+\infty$ (that is, $+\infty \not\in \TT$).
For any couple $\np{s,t} \in \bp{\RR \cup \{ +\infty \} }^2$, 
we use the notation
\begin{equation}
  \segment{s}{t}= \defset{ r \in  \TT }{ s \leq r \leq t } %\eqfinp
\end{equation}
for the \emph{segment} that joins~$s$ to~$t$
(when $s > t$, $ \segment{s}{t}= \emptyset $).

\paragraph{Special cases of discrete and continuous time.}
This setting includes the discrete time case when $\TT$ is a discrete set,
be it infinite like $\TT=\{\tinitial+k\TimeStep, \quad k \in \NN \}$
(with $\TimeStep >0$),
or finite like $\TT=\{\tinitial+k\TimeStep, \quad k =0,1,\ldots, K \}$
(with $K \in \NN$).
Of course, in the continuous time case, $\TT$ is an interval of~$\RR$,
like $[\tinitial,\horizon]$ when $\horizon <+\infty$
or $[\tinitial,+\infty[$ when $\horizon =+\infty$.
But the setting makes it possible to consider interval of 
continuous times separated by discrete times corresponding to jumps. 
For these reasons, our setting is more general than the one 
in~\cite{Rouge-Mathias-Deffuant:2013}, which considers discrete time systems.

\paragraph{Environmental illustration.}
In fisheries, investment decisions (boats, equipment) are made at large
time scale (years), regulations quotas are generally annual, 
boat operations are daily. By contrast, populations and external 
perturbations evolve in continuous time. 
Depending on the issues at hand, the modeler will choose the proper
time scales, symbolized by the set~$\TT$.
In an energy system, like a micro-grid with battery and solar panels, 
investment decisions in equipment (buying or renewal) occur at 
large time scale, whereas flows inside the system have to be decided
at short time scales (minutes).

\subsubsection{States, controls, Nature} %uncertainties}

At each time~$t \in \TT$,
\begin{itemize}
\item 
the system under consideration can be described by an 
element~$\state_{t}$ of the \emph{state set}~$\STATE_{t}$,
%equipped with a $\sigma$-field~$\tribu{\State}$.
\item 
the decision-maker (DM) makes a decision~$\control_{t}$,
taken within a \emph{control set}~$\CONTROL_{t}$.
%equipped with a $\sigma$-field~$\tribu{\Control}$.
\end{itemize}
A \emph{state of Nature}~$\omega$ affects the system, 
drawn within a \emph{sample set}~$\Omega$, also called \emph{Nature}.
No probabilistic structure is imposed on the set~$\Omega$.

\paragraph{Environmental illustration.}
In the case of dengue epidemics control at daily time steps, 
the state can be a vector of abundances 
of healthy and deseased individuals (possibly at ages),
together with the same description for the mosquito vector;
the control can be the daily fumigation effort, mosquito larva removal,
quarantine measures, or the opening and closing of sanitary facilities;
Nature represents unknown factors that affect the dengue dynamics,
like rains, humidity, mosquito biting rates, individual 
susceptibilities, etc. Some of these factirs (like rain) can be progressively 
unfolded as times passes.

\paragraph{Special case where the sample set is a set of scenarios.}
In many cases, at each time~$t \in \TT$,
an uncertainty~$\uncertain_{t}$ affects the system, 
drawn within an \emph{uncertainty set}~$\UNCERTAIN_{t}$.
Hence, a state of Nature has the form $\omega=
\sequence{\uncertain_{t}}{t \in \TT}$ --- and is called a \emph{scenario} --- 
drawn within a product {sample set} $\Omega=
\prod_{t \in \TT}\UNCERTAIN_{t}$.

\paragraph{Environmental illustration.}
The above definition of scenarios is in phase with the vocable
of scenarios in climate change mitigation; it represents sequences
of uncertainties that affect the climate evolution. 
In our framework, a scenario is not in the hands of the decision-maker;
for instance, a scenario does not include investment decisions. 

\paragraph{Relevance for resilience.}

In the case of scenarios, as the {uncertainty sets}~$\UNCERTAIN_{t}$
depend on~$t$, we cover the case where
\begin{itemize}
\item 
an uncertainty~$\uncertain_{t} \in \UNCERTAIN_{t}$ 
affects the system at each time~$t$,
possibly progressively revealed to the DM, hence available 
when he makes decisions;
\item 
other uncertainties, that are present from the start (like parameters), 
hence are part of the set~$\UNCERTAIN_{\tinitial}$;
such uncertainties are not necessarily revealed to the DM as times passes,
and remain unknown.
\end{itemize}
Our setting is more general than the one 
in~\cite{Rouge-Mathias-Deffuant:2013}.
First, we do not restrict the {sample set} to be made of scenarios
as Roug\'e, Mathias and Deffuant do.
Second, even in the case of scenarios, 
no probabilistic structure is imposed on the 
set~\( \prod_{t \in \TT}\UNCERTAIN_{t} \) 
whereas Roug\'e, Mathias and Deffuant require that it be equipped with a 
probability distribution having a density (with respect to 
an, unspecified, measure, 
likely the Lebesgue measure on a Euclidian space).

\subsubsection{State and control paths}
\label{State_and_control_paths,_scenarios_of_uncertainties,_states_of_Nature}

With the basic set~$\TT$ and the basic families of 
sets~$\sequence{\STATE_{t}}{t \in \TT}$
and $\sequence{\CONTROL_{t}}{t \in \TT}$, we define 
\begin{itemize}
\item 
the \emph{set~$\prod_{t \in \TT}\STATE_{t}$ of state paths}, made of sequences
$\sequence{\state_{t}}{t \in \TT}$ where \( \state_{t} \in \STATE_{t} \)
for all $t \in \TT$;
\emph{tail state paths} $\sequence{\state_{r}}{r \in \segment{s}{t}}$
(starting at time~$s<t$) are elements of \( \prod_{r=s}^{t}\STATE_{r} \);
\item 
the \emph{set~$\prod_{t \in \TT}\CONTROL_{t}$ of control paths}, made of sequences
$\sequence{\control_{t}}{t \in \TT}$ where \( \control_{t} \in \CONTROL_{t} \)
for all $t \in \TT$;
\emph{tail control paths} $\sequence{\control_{r}}{r \in \segment{s}{t}}$
(starting at time~$s<t$) are elements of \( \prod_{r=s}^{t}\CONTROL_{r} \).
\end{itemize}

\paragraph{Relevance for resilience.}
We introduce paths because, as stated in the abstract, our (forthcoming)
definition of resilience requires that, after any perturbation, 
the system returns to an acceptable ``regime'', that is, that 
the state-control path as a whole must return to a set of acceptable
paths (and not only the state values must belong to an acceptable subset 
of the state set). 
We introduce tail paths because resilience encapsulates the idea that recovery 
is possible after some time, and that the system remains ``normal'' after that
time.

\subsubsection{Dynamics/flow}
\label{Flow}

We now introduce a dynamics under the form of a \emph{flow} 
\( \sequence{\phi_{\segment{s}{t}}}{\np{s,t} \in \TT^2} \), that is,
a family of mappings
\begin{equation}
\phi_{\segment{s}{t}}: \STATE_{s} \times \prod_{r=s}^{t}\CONTROL_{r} 
\times \Omega 
\to \prod_{r=s}^{t}\STATE_{r} \eqfinp
\label{eq:flow} 
\end{equation}
When $s>t$, all these expressions are void because 
$ \segment{s}{t}= \emptyset $. 

The flow~$\phi_{\segment{s}{t}}$ maps 
an initial state~$\bar\state_{s} \in \STATE_{s}$ at time~$s$, 
a {tail control path} 
\( \sequence{{\control}_{r}}{r\in\segment{s}{t}} \) 
and a state of Nature~$\omega$ 
towards a tail state path
\begin{equation}
  \sequence{{\state}_{r}}{r\in \segment{s}{t}} =
\phi_{\segment{s}{t}} \bp{\bar\state_{s},
\sequence{{\control}_{r}}{r\in \segment{s}{t}},
\omega } \eqfinv
\label{eq:flow_path} 
\end{equation}
with the property that \( {\state}_{s} = \bar\state_{s} \).

\paragraph{Relevance for resilience.}
Our setting is more general than the one 
in~\cite{Rouge-Mathias-Deffuant:2013}: as illustrated below, we cover
differential and stochastic differential systems, in addition to 
iterated dynamics in  discrete time (which is the scope of 
Roug\'e, Mathias and Deffuant).
Our approach thus allows for a general treatment of resilience.

% \paragraph{Parameters in the dynamics.}
% Parameters in the dynamics can be modelled either
% as part of the mathematical expression of $\phi_{\segment{s}{t}}$,
% or as part of ${\uncertain}_{\tinitial}$.

\paragraph{Cemetery point to take into account either analytical properties
or bounds on the controls.}
In general, a state path cannot be determined by~\eqref{eq:flow_path} 
for \emph{any} state of Nature
or for \emph{any} control path, for analytical reasons (measurability,
continuity) or because of bounds on the controls.
To circumvent this difficulty, one can use a mathematical trick and 
add to any state set~$\STATE_{t}$ a cemetery point~$\partial$.
Any time a state cannot be properly defined by the flow
by~\eqref{eq:flow_path}, we attribute the 
value~$\partial$. The vocable ``cemetery'' expresses the property that,
once in the state~$\partial$, the future state values, yielded by the flow,
will all be~$\partial$. Therefore, 
the stationary state path with value~$\partial$ will be the image of those
scenarios and control paths for which no state path can be determined
by~\eqref{eq:flow_path}.

\paragraph{Special case of an iterated dynamics in  discrete time.}

In discrete time, when $\TT=\NN$, 
the flow is generally produced by the iterations of a dynamic
\begin{equation}
{\state}_{t}=\state \eqsepv {\state}_{s+1} = 
\Dynamics_{s}\np{{\state}_{s},{\control}_{s},{\uncertain}_{s}} 
\eqsepv t \geq s \eqfinp
\label{eq:iterated_dynamics_in__discrete_time}
\end{equation}

How do we include control constraints in this setting?
Suppose given a family of nonempty set-valued mappings 
\( \mathcal{\Control}_{s} : \STATE_{s} \rightrightarrows \CONTROL_{s} \), 
$s \in \TT$. 
We want to express that only controls~${\control}_{s}$ that belong
to \( \mathcal{\Control}_{s}\np{\state_{s}} \) are relevant.
For this purpose, we add to all the state sets~$\STATE_{s}$ 
a cemetery point~$\partial$. Then, when 
\( {\control}_{r} \not\in \mathcal{\Control}_{r}\np{\state_{r}} \)
in~\eqref{eq:iterated_dynamics_in__discrete_time}
for at least one~\( r \in \segment{s}{t} \), we set
\( \phi_{\segment{s}{t}} \bp{\bar\state_{s},
\sequence{{\control}_{r}}{r\in \segment{s}{t}},
\sequence{{\uncertain}_{r}}{r\in \segment{s}{t}},\gamma } = 
\sequence{\partial}{s\in \segment{t}{\horizon}} \) 
in~\eqref{eq:flow_path}.

\paragraph{Environmental illustration.}
In natural resource management, many population models 
(anmal, plants) are given by discrete time abundance-at-age 
dynamical equations.
Outside population models, many stock problems are also based upon
discrete time dynamical equations. This is the case of dam management,
where water stock balance equations are written at a daily scale
(possibly less like every eight hours, or possibly more  like
months for long term planning);
control constraints represent the properties that turbined water must be 
less than the current water stock and bounded by turbine capacity.
%The time step is generally the one of the decision process.

\paragraph{Special case of differential systems.}

In continuous time, the mapping $ \phi_{\segment{s}{t}}$ in~\eqref{eq:flow}
generally cannot be defined over the whole set 
\( \prod_{r=s}^{t}\CONTROL_{r} \times \Omega \). 
Tail control paths and states of Nature need to be restricted 
to subsets of $\prod_{r=s}^{t}\CONTROL_{r}$ and $\Omega$,
like the continuous ones for example when dealing with Euclidian spaces. 
For instance, when $\TT=\RR_+$ and the flow is produced by a smooth 
dynamical system on a Euclidian space~$\STATE=\RR^n$
\begin{equation}
{\state}_{t}=\state \eqsepv 
\dot{{\state}}_{s} =\dynamics_{s}\np{{\state}_{s},{\control}_{s}}
\eqsepv s \geq t \eqfinv
\end{equation}
control paths  \( \sequence{{\control}_{s}}{s\in\segment{t}{\horizon}} \) 
are generally restricted to piecewise continuous ones
for a solution to exist. 

\paragraph{Special case of stochastic differential equations.}

Under certain technical assumptions, 
a stochastic differential equation
\begin{equation}
  d{\va{\State}}_{s} =\dynamics_{s}
\np{\va{\State}_{s},\va{\Control}_{s},\va{\Uncertain}_{s}}ds +
g\np{\va{\State}_{s},\va{\Control}_{s},\va{\Uncertain}_{s}} d\va{\Uncertain}_{s}
\eqfinv
\end{equation}
where \( \sequence{\va{\Uncertain}_{s}}{s\in \RR_+} \) is a Brownian motion,
gives rise to solutions in the strong sense. In that case, 
a flow can be defined (but not over the whole set 
\( \prod_{r=s}^{t}\CONTROL_{r} \times \Omega \)).

\paragraph{The case of the history flow.}
Any possible state derives from the so-called \emph{history}
\begin{equation}
  \history_t=\bp{\sequence{{\control}_{r}}{r\in \segment{\tinitial}{t}},
\omega } \eqfinp
\end{equation}
In that case, the flow~\eqref{eq:flow_path} is trivially given by 
\(  \sequence{{\history}_{r}}{r\in \segment{s}{t}} =
\bp{\history_{s},
\sequence{{\control}_{r}}{r\in \segment{s}{t}}  } \).
We will use the notion of history when we compare our approach with the 
viability approach to resilience.

\subsection{Adapted and admissible strategies}

A control~${\control}_{t}$ is an element of the {control set}~$\CONTROL_{t}$.
A \emph{policy} (at time~$t$) is a mapping 
\begin{equation}
  \strategy_{t} : \STATE_{t} \times \Omega \to \CONTROL_{t} 
\label{eq:policy}
\end{equation}
with image in the {control set}~$\CONTROL_{t}$.
A \emph{strategy} is a sequence 
\( \sequence{\strategy_{t}}{t\in \TT} \) of policies. 

\paragraph{Environmental illustration.}
In climate change mitigation, a strategy can be an investment policy 
in renewable energies as a function of the past observed temperatures.
In epidemics control, a strategy can be quarantine measures or
vector control as a function of observed infected individuals.

\subsubsection{Admissible strategies} 

Suppose given a family of nonempty set-valued mappings 
\( \mathcal{\Control}_{t} : \STATE_{t} \times \Omega 
\rightrightarrows \CONTROL_{t} \), $t \in \TT$. 
An \emph{admissible strategy} is a strategy
\( \sequence{\strategy_{t}}{t\in \segment{\tinitial}{\horizon}} \) 
such that control constraints are satisfied in the sense that,
for all $t \in \TT$, 
\begin{equation}
\strategy_{t}\np{\state_{t}, \omega }
\in \mathcal{\Control}_{t}\np{\state_{t},\omega } %
\eqsepv \forall 
\np{ \state_{t}, \omega }
\in \STATE_{t} \times \Omega \eqfinp
\end{equation}

\subsubsection{Adapted strategies}

Suppose that the sample set~$\Omega$ is equipped with a 
\emph{filtration}
\( \sequence{\tribu{F}_{t}}{t\in \TT} \).
Hence each $\tribu{F}_{t}$ is a $\sigma$-field and 
the sequence \( t \mapsto \tribu{F}_{t} \) is increasing (for the inclusion
order). 
Suppose that each \emph{state set}~$\STATE_{t}$,
is equipped with a $\sigma$-field~$\tribu{\State}_{t}$.

An \emph{adapted policy} is a mapping~\eqref{eq:policy}
which is measurable with respect to the product
$\sigma$-field~$\tribu{\State}_{t} \otimes \tribu{F}_{t} $.
An \emph{adapted strategy} is a family 
\( \sequence{\strategy_{t}}{t\in \segment{\tinitial}{\horizon}} \) 
of adapted policies.

\paragraph{Special case where the sample set is a set of scenarios.}
Consider the case where $\Omega= \prod_{t \in \TT}\UNCERTAIN_{t}$ and 
where each set~$\UNCERTAIN_{t}$ is equipped with a $\sigma$-field 
$\tribu{\Uncertain}_{t}$ (supposed to contain the singletons).
The natural {filtration} \( \sequence{\tribu{F}_{t}}{t\in \TT} \)
is given by 
\begin{equation}
  \tribu{F}_{t} = \bigotimes_{r \leq t} \tribu{\Uncertain}_{r}
\otimes \bigotimes_{s > t} \{ \emptyset, \UNCERTAIN_{s} \} \eqfinp 
\end{equation}
Then, in an {adapted strategy} 
\( \sequence{\strategy_{t}}{t\in \segment{\tinitial}{\horizon}} \),
each policy can be identified with a mapping of the form
\cite{Carpentier-Chancelier-Cohen-DeLara:2015}
\begin{equation}
  \strategy_{t} : \STATE_{t} \times \prod_{r=\tinitial}^{t}\UNCERTAIN_{r} 
\to \CONTROL_{t} \eqfinp
\end{equation}
In that case, 
our definition of adapted strategy means that the DM can, at time~$t$, 
use no more than 
time~$t$, current state value~$\state_{t}$ and past scenario 
$\sequence{{\uncertain}_{s}}{s\in \segment{\tinitial}{t} }$ 
to make his decision \( \control_{t} = \strategy_{t}\np{ \state_{t},
\sequence{ {\uncertain}_{s}}{s\in \segment{\tinitial}{t} } } \).

\paragraph{Relevance for resilience.}
Though this is not the most general framework to handle information
(see \cite{Carpentier-Chancelier-Cohen-DeLara:2015} for a more general
treatment of information),
we hope it can enlighten the notion of \emph{adaptive response} often
found in the resilience literature.

Our setting is more general than the one 
in~\cite{Rouge-Mathias-Deffuant:2013}: indeed, 
Roug\'e, Mathias and Deffuant only consider state feedbacks,
that is, Markovian strategies as defined below.
By contrast, our setting includes the case of corrupted and 
partially observed state feedback strategies, that is, the case where 
strategies have as input a partial observation of the state that is
corrupted by noise.

\paragraph{Special case of Markovian or state feedback strategies.}

Markovian or state feedback policies are of the form 
\begin{equation}
  \strategy_{t} : \STATE_{t} \to \CONTROL_{t} \eqfinp
\end{equation}
With this definition, we express that, at time~$t$, the DM can only use
time~$t$ and current state value~$\state_{t}$ --- but not the 
state of Nature~$\omega$ ---
to make his decision \( \control_{t} = \strategy_{t}\np{ \state_{t} } \).
In some cases (when dynamic programming applies for instance), it is enough to
restrict to Markovian strategies, much more economical than general strategies.

\subsection{Closed loop flow} 

From now on, when we say ``strategy'', we mean ``adapted and admissible
strategy''.
\bigskip

Given an initial state and a state of Nature, 
a strategy will induce a state path
thanks to the flow: this gives the closed loop flow as follows. 

Let $s \in \TT$ and $t \in \TT$, with $s < t$.
Let \( \sequence{\strategy_{t}}{t\in \TT} \) be a strategy. 
We suppose that, for any initial state~$\bar\state_{s} \in \STATE_{s}$
and any state of Nature~$\omega$, the following system of 
(closed loop) equations
\begin{subequations}
  \begin{align}
      \sequence{{\state}_{r}}{r\in \segment{s}{t}} &=
\phi_{\segment{s}{t}} \bp{\bar\state_{s},
\sequence{{\control}_{r}}{r\in \segment{s}{t}}, \omega } \\
{\control}_{r} &= \strategy_{r}\np{ {\state}_{r} , \omega }
\eqsepv \forall r \in \segment{s}{t} 
  \end{align}
\end{subequations}
has a \emph{unique} solution 
\( \bp{ \sequence{{\state}_{r}}{r \in \segment{s}{t}} , 
\sequence{{\control}_{r}}{r \in \segment{s}{t}} } \).
Quite naturally, we define the 
\emph{closed loop flow}~\( \phi_{\segment{s}{t}}^{\strategy} \) by 
\begin{equation}
  \phi_{\segment{s}{t}}^{\strategy}
  \bp{\bar\state_{s},\omega }
= \bp{ \sequence{{\state}_{s}}{s\in \segment{s}{t}} , 
\sequence{{\control}_{s}}{s\in \segment{s}{t}} } \eqfinp
\label{eq:closed_loop_flow}
\end{equation}
% For flows given by dynamical systems in continuous or discrete time, uniqueness
% is easily obtained. 

\section{Resilience: a mathematical framework}
\label{Resilience:_A_Mathematical_Framework}

Equipped with the material in~Sect.~\ref{Ingredients_for_an_abstract_control_system_with_uncertainties},
we now frame the concept of resilience in mathematical garbs.
For this purpose, we introduce the notion of \emph{recovery regime}.
Compared to other definitions of resilience 
\cite{Holling:73,Martin:2004,Martin-Deffuant-Calabrese:2011,Rouge-Mathias-Deffuant:2013},
our definition requires that, after any perturbation, 
the state-control path as a whole can be driven, by a proper strategy, 
to a set of acceptable paths 
(and not only the state values must belong to an acceptable subset 
of the state set, asymptotically or not). 
In addition, as state and control paths are contingent on uncertainties, 
we require that they lay within an acceptable domain of random processes, 
called recovery regimes.

Once again, as the reference \cite{Rouge-Mathias-Deffuant:2013} is the more
mathematically driven paper on resilience, we will systematically 
emphasize in what our approach differs from that of 
Roug\'e, Mathias and Deffuant.

\subsection{Robustness, resilience and random processes}

The notion of robustness captures a form of stability to perturbations;
it is a static notion, as no explicit reference to time is required.
By contrast, the concept of resilience makes reference to 
time (dynamics), strategies and perturbations. 
This is why, to speak of resilience --- 
a notion that mixes time and randomness ---
we find it convenient to use the framework of 
random processes, although this does not mean that we require any
probability.

From now on, we consider that the {sample space}~$\Omega$ is a
measurable set equipped with a $\sigma$-field~$\tribu{F}$ 
(but not necessarily equipped with a probability).
When we consider a deterministic setting, $\Omega$ is reduced to a singleton
(and ignored).
The set of \emph{measurable mappings} from~$\Omega$ to 
any measurable set~$\YY$ will be denoted by \( \Adapted{\Omega}{\YY} \).
Elements of \( \Adapted{\Omega}{\YY} \) are called \emph{random variables}
or \emph{random processes}, although this does not imply the existence of
an underlying probability.
Random variables are designated with bold capital letters like~$\va{Z}$.
From now on, every {state set}~$\STATE_{t}$ is a measurable set 
equipped with a $\sigma$-field~$\tribu{\State}_{t}$,
every {control set}~$\CONTROL_{t}$ with a $\sigma$-field~$\tribu{\Control}_{t}$,
and, when needed, every {uncertainty set}~$\UNCERTAIN_{t}$ 
with a $\sigma$-field $\tribu{\Uncertain}_{t}$.

Fields are introduced when probabilities are needed.
When they are not, as with the robust setting, it suffices to equip 
all sets with their complete $\sigma$-fields, made of all subsets.
Then, {measurable mappings} \( \Adapted{\Omega}{\YY} \)
from~$\Omega$ to any set~$\YY$ are all mappings.

% \paragraph{Example of shocks.}

% Suppose that the {sample space}~$\Omega$ is a measurable set
% (equipped with a $\sigma$-field~$\tribu{F}$) 
% with a probability~$\PP$.
% Suppose that $\TT=\RR_+$, $\UNCERTAIN_{t}=\RR$, for all $t \in \TT$,
% and that $\Gamma$ is reduced to a singleton, hence ignored. 
% Let $\TimeStep >0$ represent a small time period, 
% $\bar\Uncertain$ a high value for uncertainties,
% and $\beta < 1 $, $\beta \approx 1 $ a probability level. 
% The set 
% \begin{align}
%   \Shock^{t} = \{ \sequence{\va{\Uncertain}_{s}}{s \geq t } 
% \in \Adapted{\Omega}{\prod_{r=t}^{\horizon}\UNCERTAIN_{r}} \mid 
% % \in \Adapted{\segment{t}{\horizon}}{\UNCERTAIN} \mid 
% & \va{\Uncertain}_{s}= 0  \eqsepv \forall s > t+\TimeStep \eqsepv \nonumber \\
% & \PP\bc{\va{\Uncertain}_{s} \geq
% \bar\Uncertain  \eqsepv \forall s \in [t, t+\TimeStep ] } \geq \beta \} 
% \end{align}
% % \begin{equation}
% %   \Shock^{t} = \defset{ \sequence{\va{\Uncertain}_{s}}{s \geq t } }{%
% % \va{\Uncertain}_{s}= 0  \eqsepv \forall s > t+\TimeStep \eqsepv \PP\bp{\va{\Uncertain}_{s} \geq
% % \bar\Uncertain  \eqsepv \forall s \in [t, t+\TimeStep ] } \geq \beta } 
% % \end{equation}
% contains stochastic process that, with high probability at least~$\beta$,
% display high values at least~$\bar\Uncertain$ for times close to~$t$.

\subsection{Recovery regimes}

\emph{Recovery regimes}, starting from $t \in \TT$,
are subsets of random processes of the form
\begin{equation}
  \Acceptable^{t} \subset 
\Adapted{\Omega}{\prod_{s=t}^{\horizon}\STATE_{s} \times 
\prod_{s=t}^{\horizon}\CONTROL_{s} }
% \Adapted{\segment{t}{\horizon}}{\STATE\times\CONTROL} 
\eqfinp
\end{equation}
When there are no uncertainties, $\Omega$ is reduced to a singleton, so that 
\( \Acceptable^{t} \subset 
\prod_{s=t}^{\horizon}\STATE_{s} \times \prod_{s=t}^{\horizon}\CONTROL_{s} \),
as in the two first following examples.

\paragraph{Example of recovery regimes converging to an equilibrium.}

Let $\TT=\RR_+$, $\STATE_{t}=\RR^n$ and $\CONTROL_{t}=\RR^m$, 
for all $t \in \TT$. 
Let \( \bar\state \) be an equilibrium point of the dynamical system
\(  \dot{\state}_{s} =\dynamics\np{\state_{s},\bar\control} \) 
when the control is stationary equal to~$\bar\control$, that is,
\( 0 =\dynamics\np{\bar\state,\bar\control} \).
The recovery regimes, starting from $t \in \TT$,
converging to the equilibrium~$\bar\state$ 
form the set
\begin{equation}
  \Acceptable^{t} = \defset{ \bp{ \sequence{\state_{s}}{s \geq t},
\sequence{\control_{s}}{s \geq t} } \in \STATE^{[t,+\infty[} \times
\CONTROL^{[t,+\infty[} }%
{\state_{s} \to_{s\to +\infty} \bar\state } \eqfinp
\end{equation}
In general, the equilibrium~$\bar\state$ is supposed to be 
asymptotically stable, locally or globally.

A more general definition would be 
\begin{equation}
  \Acceptable^{t} = \defset{ \bp{ \sequence{\state_{s}}{s \geq t},
\sequence{\control_{s}}{s \geq t} } }%
{\lim_{s\to +\infty}\state_{s} \textrm{ exists} } \eqfinv
\end{equation}
and, to account for constraints on the values taken by the controls, 
we can consider
\begin{equation}
  \Acceptable^{t} = \defset{ \bp{ \sequence{\state_{s}}{s \geq t},
\sequence{\control_{s}}{s \geq t} } }%
{\lim_{s\to +\infty}\state_{s} \textrm{ exists and } 
\control_{s} \in \mathcal{\Control}_{s}\np{\state_{s}} \eqsepv \forall s \geq t } \eqfinv
\end{equation}
where \( \mathcal{\Control}_{s} : \STATE_{s} \rightrightarrows \CONTROL_{s} \), 
for all $s \in \TT$.

\paragraph{Example of bounded recovery regimes.}

Let $\TT=\RR_+$, $\STATE_{t}=\RR^n$ and $\CONTROL_{t}=\RR^m$, 
for all $t \in \TT$. 
If~$B$ is a bounded region of $\STATE=\RR^n$, we consider 
\begin{equation}
 \Acceptable^{t} = \defset{ \bp{ \sequence{\state_{s}}{s \geq t},
\sequence{\control_{s}}{s \geq t} } }%
{\state_{s} \in B \eqsepv \forall s \geq t } \eqfinp
\end{equation}
When $B$ is a ball of small radius~$\rho>0$ around
the equilibrium~$\bar\state$, we obtain state paths that remain 
close to~$\bar\state$.

\paragraph{Example of random recovery regimes.}

Suppose that the measurable {sample space}~$(\Omega,\tribu{F})$ is 
equipped with a probability~$\PP$. 
Let $\STATE_{t}=\RR^n$ and $\CONTROL_{t}=\RR^m$, 
for all $t \in \TT$. 
Letting~$B$ be a bounded region of $\STATE=\RR^n$ and 
$\beta \in ]0,1[ $, the set 
\begin{align}
  \Acceptable^{t} = & \{ \bp{ \sequence{\va{\State}_{s}}{s \geq t },
\sequence{\va{\Control}_{s}}{s \geq t} } 
\in \Adapted{\Omega}{\prod_{s=t}^{\horizon}\STATE_{s} \times 
\prod_{s=t}^{\horizon}\CONTROL_{s} } \mid  \nonumber \\
& \PP\bc{ \exists s \geq t \mid \va{\State}_{s} \not \in B } \leq \beta \} 
\end{align}
represents state paths that get at least once outside the bounded region~$B$
with a probability less than~$\beta$. 

If $\TT$ is discrete, the set 
\begin{align}
  \Acceptable^{t} = & \{ \bp{ \sequence{\va{\State}_{s}}{s \geq t },
\sequence{\va{\Control}_{s}}{s \geq t} } 
\in \Adapted{\Omega}{\prod_{s=t}^{\horizon}\STATE_{s} \times 
\prod_{s=t}^{\horizon}\CONTROL_{s} } \mid  \nonumber \\
& \PP\bc{ \exists s_1 \geq t \eqsepv s_2 \geq t \eqsepv s_3 \geq t 
\mid \va{\State}_{s_1} \not \in B \eqsepv
\va{\State}_{s_2} \not \in B \eqsepv \va{\State}_{s_3} \not \in B } = 0 \} 
\end{align}
represents state paths that get no more than two times outside the bounded region~$B$.

\subsection{Resilient strategies and resilient states}
\label{Resilient_strategies_and_resilient_states}

Consider a starting time~$t\in\TT$ 
and an initial state~$\bar\state_{t}\in\STATE_{t}$. 
We say that the strategy~\( \strategy \)
is a \emph{resilient strategy} starting from time~$t$ 
in state~$\bar\state_{t}$ if the random process 
\( \Bp{ \sequence{\va{\State}_{s}}{s\in \segment{t}{\horizon}} ,
\sequence{\va{\Control}_{s}}{s\in \segment{t}{\horizon}}  } \)
given by 
\begin{subequations}
  \begin{align}
      \sequence{\va{\State}_{s}\np{\omega}}{s\in \segment{t}{\horizon}} &=
  \phi_{\segment{t}{\horizon}}^{\strategy}\bp{\state_{t},\omega }
\\
\va{\Control}_{s}\np{\omega} &= \strategy_{s}\np{ \va{\State}_{s}\np{\omega} , 
\omega }
\eqsepv \forall s\in \segment{t}{\horizon} \eqfinv
  \end{align}
\label{eq:output}
\end{subequations}
where the closed loop flow~$\phi_{\segment{t}{\horizon}}^{\strategy}$
is given in~\eqref{eq:closed_loop_flow},
is such that
\begin{equation}
\Bp{ \sequence{\va{\State}_{s}}{s\in \segment{t}{\horizon}} ,
\sequence{\va{\Control}_{s}}{s\in \segment{t}{\horizon}}  } 
\in \Acceptable^{t} \eqfinp  
\end{equation}
Notice that we do not use the part 
\( \sequence{\strategy_{r}}{r < t} \) of
the strategy~\( \strategy = 
\sequence{\strategy_{r}}{r\in \segment{\tinitial}{\horizon} } \).

With this definition, a resilient strategy is able to drive the 
state-control random process into an acceptable regime.
As a resilient strategy is adapted, it can ``adapt'' to the 
past values of the randomness but no to its future values
(hence, our notion of resilience does not require clairvoyance of the DM). 
Our definition of resilience 
is a form of controlability for whole random processes (regimes):
a resilient strategy has the property 
to shape the closed loop flow~$\phi_{\segment{t}{\horizon}}^{\strategy}$
so that it belongs to a given subset of random processes.

We denote by~\( \ResilientStrategies_{t}\np{\bar\state_{t}} \)
the set of \emph{resilient strategies} at time~$t$,
starting from state~$\bar\state_{t}$.
The set of \emph{resilient states} at time~$t$ is
\begin{equation}
  \ResilientStates_{t} = \defset{ \bar\state_{t} \in \STATE_{t} }%
{\ResilientStrategies_{t}\np{\bar\state_{t}} \not = \emptyset } \eqfinp 
\end{equation}

\section{Illustrations}
\label{Illustrations}

In~Sect.~\ref{Resilience:_A_Mathematical_Framework}, 
we have provided some illustrations in the course of the exposition. 
Now, we make the connection between the previous setting and
two other settings, the resilience ``\`a la Holling'' \cite{Holling:73} 
in~\S\ref{Holling}
and the resilience-viability framework
\cite{Martin:2004,Martin-Deffuant-Calabrese:2011,Rouge-Mathias-Deffuant:2013} 
in~\S\ref{resilience-viability}.

\subsection{Deterministic control dynamical system with attractor}
\label{Holling}

As the paper~\cite{Holling:73} does not contain a single equation, 
it is bit risky to force the seminal Holling's contribution into our setting.
However, it is likely that it corresponds to $\TT=\RR_+$ and 
to recovery regimes of the form 
\begin{equation}
  \Acceptable^{t} = \defset{ \bp{ \sequence{\state_{s}}{s \geq t},
\sequence{\control_{s}}{s \geq t} } }%
{\state_{s} \textrm{ converges towards an attractor} } \eqfinp
\end{equation}
Note that, as often in the ecological literature on resilience
\cite{Arnoldi-Loreau-Haegeman:2016},
the underlying dynamical system is not controlled.

\subsection{Resilience and viability}
\label{resilience-viability}

Some authors 
\cite{Martin:2004,Martin-Deffuant-Calabrese:2011,Rouge-Mathias-Deffuant:2013} 
propose to frame resilience within the mathematical theory of
viability \cite{Aubin:1991}.

Let $\STATE_{t}=\STATE$ and $\CONTROL_{t}=\CONTROL$, for all $t \in \TT$. 
Let $\sdo \subset \STATE$ denote a set made of ``acceptable states''.
Let \( \mathcal{\Control}_{s} : \STATE \rightrightarrows \CONTROL \), 
$s \in \TT$
be a family of set-valued mappings that represent control constraints.

\subsubsection{Deterministic viability}

Consider a starting time~$t\in\TT$ and the recovery regimes
\begin{align}
 \Acceptable^{t} = \{ \bp{ \sequence{\state_{s}}{s \geq t},
\sequence{\control_{s}}{s \geq t} } \mid 
\state_{s} \in \sdo \eqsepv  \control_{s} \in \mathcal{\Control}_{s}\np{\state_{s}} 
\eqsepv \forall s \geq t \}  \eqfinp
\end{align}
Then, a resilient strategy is one that is able to drive the 
state towards the set~$\sdo$ of acceptable states.

\subsubsection{Robust viability}
\label{Deterministic_and_robust_viability}

When there are no uncertainties, we just established 
a connection between recovery regimes and viability.
But, with uncertainties, as resilience requires a form of 
stability ``whatever the perturbations'', 
we are in the realm of \emph{robust viability} \cite{DeLara-Doyen:2008},
as follows. 

Let  \( \overline\Omega \subset \Omega \), corresponding to the (nonempty) subset of 
states of Nature with respect to which the DM expects the system to be robust.
Consider a starting time~$t\in\TT$ and the recovery regimes
\begin{align}
 \Acceptable^{t} = \{ & 
\bp{ \sequence{\va{\State}_{s}}{s\in \segment{t}{\horizon}}
\sequence{\va{\Control}_{s}}{s\in \segment{t}{\horizon}} } 
\in \Adapted{\Omega}{\prod_{s=t}^{\horizon}\STATE_{s} \times 
\prod_{s=t}^{\horizon}\CONTROL_{s} } \mid \nonumber \\
& \exists \va{\tau} \in \Adapted{\Omega}{\TT} \eqsepv 
\forall \omega \in \overline\Omega \eqsepv \nonumber \\
& \va{\tau}\np{\omega} \geq t  \eqsepv 
\va{\State}_{s}\np{\omega} \in \sdo \eqsepv 
\va{\Control}_{s}\np{\omega} \in 
\mathcal{\Control}_{s}\np{\va{\State}_{s}\np{\omega}} 
\eqsepv \forall s \geq \va{\tau}\np{\omega} \} \eqfinp
\label{eq:recovery_regimes_robust_viability}
\end{align}
Then, a resilient strategy is one that is able to drive the 
state towards the set~$\sdo$ of acceptable states, 
after a random time~$\va{\tau}$, 
whatever the perturbations in~$\overline\Omega$.

\subsubsection{Robust viability and recovery time attached to a resilient strategy }
\label{recovery_time}

Let  \( \overline\Omega \subset \Omega \), whose elements 
can be interpreted as shocks.
Consider a starting time~$t\in\TT$ 
and an initial state~$\bar\state_{t}\in\STATE_{t}$. 
If \( \strategy=\sequence{\strategy_{s}}{s\in \segment{t}{\horizon}} \)
is a resilient strategy for the 
recovery regimes~\eqref{eq:recovery_regimes_robust_viability},
the \emph{recovery time} is the random time defined by 
\begin{align}
 \va{\tau}\np{\omega} =\inf\{ r \in \segment{t}{\horizon} \mid &
\sequence{\va{\State}_{s}\np{\omega}}{s\in \segment{t}{\horizon}} =
\phi_{\segment{t}{\horizon}}^{\strategy}  \bp{\bar\state_{t}, \omega} \nonumber \\
& \va{\Control}_{s}\np{\omega} = 
\strategy_{s}\np{ \va{\State}_{s}\np{\omega} , \omega } 
\eqsepv \forall s\in \segment{t}{\horizon}  \nonumber \\
& \va{\State}_{s}\np{\omega} \in \sdo \eqsepv 
\va{\Control}_{s}\np{\omega} \in 
\mathcal{\Control}_{s}\np{ \va{\State}_{s}\np{\omega} }
\eqsepv \forall s \geq r  \} \eqfinv
\end{align}
for all \( \omega \in \Omega \), 
with the convention that $\inf \emptyset = +\infty$. 

Thus, the resilient strategy drives the 
state towards the set~$\sdo$ of acceptable states, 
after the random time~$\va{\tau}$, 
whatever the perturbations (shocks) in~$\overline\Omega$. 
By contrast, the so-called \emph{time of crisis}
occurs before~$\va{\tau}$ \cite{Doyen-Saint-Pierre:1997}.

\subsubsection{Stochastic viability}
\label{Stochastic_viability}

Suppose that the measurable {sample space}~$(\Omega,\tribu{F})$ is 
equipped with a probability~$\PP$
and let $\beta \in [0,1] $, represent a probability level. 
Consider a starting time~$t\in\TT$ and the recovery regimes
\begin{align}
  \Acceptable^{t} = \{ & 
\bp{ \sequence{\va{\State}_{s}}{s\in \segment{t}{\horizon}},
\sequence{\va{\Control}_{s}}{s\in \segment{t}{\horizon}} }
\in \Adapted{\Omega}{\prod_{s=t}^{\horizon}\STATE_{s} \times 
\prod_{s=t}^{\horizon}\CONTROL_{s} } \mid  \nonumber \\
& \PP\bc{ \va{\State}_{s} \in \sdo \eqsepv 
\va{\Control}_{s} \in \mathcal{\Control}_{s}\np{ \va{\State}_{s} }
\eqsepv \forall s \geq t } \geq \beta \} \eqfinp
\end{align}
With these recovery regimes, 
we express that the probability to satisfy state and control constraints after
time~$t$ is at least~$\beta$ \cite{Doyen-DeLara:2010}.

\subsubsection{Discussion and comparison with the viability theory approach for resilience.}

Our setting is more general than the viability theory approach for resilience
as introduced in 
\cite{Martin:2004,Martin-Deffuant-Calabrese:2011,Rouge-Mathias-Deffuant:2013}.
Indeed, the viability approach to resilience deals with constraints time 
by time; our approach does not. 

To illustrate our point, consider
the deterministic case with discrete and finite time,
and scalar controls, to make things easy.
It is clear that the recovery regimes given by
\begin{subequations}
\begin{align}
 \Acceptable^{t} &= \{ \bp{ \sequence{\state_{s}}{s \geq t},
\sequence{\control_{s}}{s \geq t} } \mid 
\min_{s \geq t} \control_{s} \leq 0 \} \\
&= \{ \bp{ \sequence{\state_{s}}{s \geq t},
\sequence{\control_{s}}{s \geq t} } \mid 
\exists s \geq t \eqfinv \control_{s} \leq 0 \}  \eqfinp
\end{align}  
\end{subequations}
cannot be expressed as time by time constraints on the controls. 

Of course, the viability approach \emph{could} handle such a case,
but at the price of extending the state and the dynamics, to turn an 
intertemporal constraint into a time by time constraint.
For instance, with the history state introduced at the end of~\S\ref{Flow},
we can always express any recovery regimes set 
as viability constraints.
In the example above, we do not need the whole history to turn
the set~$\Acceptable^{t}$ into one described by time by time constraints:
it suffices to introduce an additional component to the state like
\( \sum_{s \geq t} \1_{ \{ \control_{s} \leq 0 \} } \) in discrete time
%(or an integral in continuous time), 
and impose the final constraint that this new part of an extended state
be non zero.

To sum up, our approach to resilience covers more recovery regimes,
described with the original states and controls,
than those captured by the time by time constraints that make the 
specificity of the viability approach to resilience.

\section{Resilience and risk}
\label{Resilience_and_risk_control/minimization}

We now sketch how concepts from risk measures (introduced initially in
mathematical finance \cite{Follmer-Schied:2002}) can be imported to tackle
resilience issues.
This again is a novelty with respect to the 
stochastic viability theory approach 
for resilience as in~\cite{Rouge-Mathias-Deffuant:2013}.
Risk measures are potential candidates as \emph{indicators of resilience}.

\subsection{Recovery regimes given by risk measures}

We start by a definition of recovery regimes given by risk measures,
then we provide examples. 

\subsubsection{Definition of recovery regimes given by extended risk measures}

Suppose that \( \TT \subset \RR \) is equipped with the trace~$\tribu{T}$ 
of the Borel field of~$\RR$. Then, \( \TT \times \Omega \) 
is a measurable space when equipped with the product $\sigma$-field
$\tribu{T} \otimes \tribu{F}$.
Then, any random process in 
\( \Adapted{\Omega}{\prod_{s=t}^{\horizon}\STATE_{s} \times 
\prod_{s=t}^{\horizon}\CONTROL_{s} } \)
can be identified with a random variable in 
\( \Adapted{\segment{t}{\horizon} \times 
\Omega}{\bigcup_{s=t}^{\horizon}\STATE_{s} \bigcup
\bigcup_{s=t}^{\horizon}\CONTROL_{s} } \).
We call \emph{extended risk measure} any 
$\GG_t$ that maps random variables 
in \( \Adapted{\segment{t}{\horizon} \times \Omega}%
{\bigcup_{s=t}^{\horizon}\STATE_{s} \bigcup
\bigcup_{s=t}^{\horizon}\CONTROL_{s} } \)
towards the real numbers \cite{Follmer-Schied:2002}.
The lower the risk measure~$\GG_t$, the better.

The basic example of a risk measure is the mathematical expectation 
under a given probability distribution. 
A celebrated risk measure in mathematical finance is the 
\emph{tail/average/conditional value-at-risk}.

With $\GG_t$ an extended risk measure and 
$\alpha \in \RR$ a given \emph{risk level},
we define recovery regimes by
\begin{align}
  \Acceptable^{t} = \{ & 
\bp{ \sequence{\va{\State}_{s}}{s\in \segment{t}{\horizon}},
\sequence{\va{\Control}_{s}}{s\in \segment{t}{\horizon}} }
\in  \Adapted{\Omega}{\prod_{s=t}^{\horizon}\STATE_{s} \times 
\prod_{s=t}^{\horizon}\CONTROL_{s} } \mid \nonumber \\
& \GG_{t} \bc{ \sequence{\va{\State}_{s}}{s\in \segment{t}{\horizon}},
\sequence{\va{\Control}_{s}}{s\in \segment{t}{\horizon}} } 
\leq \alpha \} \eqfinp 
\label{eq:alpha}
\end{align}
The quantity 
\( \GG_{t} \bc{ \sequence{\va{\State}_{s}}{s\in \segment{t}{\horizon}},
\sequence{\va{\Control}_{s}}{s\in \segment{t}{\horizon}} } \) measures 
the ``risk'' borne by the random process
\( \bp{ \sequence{\va{\State}_{s}}{s\in \segment{t}{\horizon}},
\sequence{\va{\Control}_{s}}{s\in \segment{t}{\horizon}} } \).
Therfore, recovery regimes like in~\eqref{eq:alpha}
represent a form of ``risk containment'' under the level~$\alpha$.

\subsubsection{Robust viability and the worst case risk measure}

The robust viability inspired definition of resilience 
in~\S\ref{Deterministic_and_robust_viability} corresponds to~\eqref{eq:alpha}
with \( \alpha < 1 \) and the \emph{worst case risk measure}
\begin{equation}
 \GG_{s}\bp{ \sequence{\va{\State}_{s}}{s\in \segment{t}{\horizon}},
\sequence{\va{\Control}_{s}}{s\in \segment{t}{\horizon}} } = 
\sup_{s\in \segment{t}{\horizon}} \sup_{\omega \in \overline\Omega} 
\1_{\sdo^c}\bp{\va{\State}_{s}\np{\omega}} \eqfinv
\end{equation}
where \( \overline\Omega \subset \Omega \).
Indeed, \( \GG_{t} \bc{ \sequence{\va{\State}_{s}}{s\in \segment{t}{\horizon}},
\sequence{\va{\Control}_{s}}{s\in \segment{t}{\horizon}} } 
\leq \alpha < 1 \) means that \( \1_{\sdo^c}\bp{\va{\State}_{s}\np{\omega}}
\equiv 0 \), that is, the state \( \bp{\va{\State}_{s}\np{\omega}} \)
always belongs to~$\sdo$ (as $\sdo^c$ is the complementary set of~$\sdo$
in~$\STATE$) for all $\omega \in \overline\Omega$.

Here, the {worst case risk measure} captures that
the state \( \va{\State}_{s}\np{\omega} \) belongs to~$\sdo$
both for all times --- the core of viability, here handled by the 
term $\sup_{s \geq t} $ --- 
and for all states of Nature in~$\overline\Omega$ --- the core of robustness, here handled by the 
term $\sup_{\omega \in \overline\Omega}$. 

\subsubsection{Stochastic viability and beyond: ambiguity}

The stochastic viability inspired definition of resilience 
in~\S\ref{Stochastic_viability} corresponds to~\eqref{eq:alpha}
with \( \alpha=1-\beta \) and the risk measure
\begin{equation}
 \GG_{s}\bp{ \sequence{\va{\State}_{s}}{s\in \segment{t}{\horizon}},
\sequence{\va{\Control}_{s}}{s\in \segment{t}{\horizon}} } = 
\PP\bc{ \exists s \geq t \mid \va{\State}_{s} \not\in \sdo \mtext{ or }
\va{\Control}_{s} \not\in \mathcal{\Control}_{s}\np{ \va{\State}_{s} } } 
\eqfinp 
\end{equation}
Now, suppose that different risk-holders do not share the same beliefs
and let $\mathcal{P}$ denote a set of probabilities on 
$\np{\Omega, \tribu{F}}$.
We can arrive at an \emph{ambiguity viability} inspired definition 
of resilience using the risk measure
\begin{equation}
 \GG_{s}\bp{ \sequence{\va{\State}_{s}}{s\in \segment{t}{\horizon}},
\sequence{\va{\Control}_{s}}{s\in \segment{t}{\horizon}} } = \sup_{\PP\in \mathcal{P}}
\PP\bc{ \exists s \geq t \mid \va{\State}_{s} \not\in \sdo \mtext{ or }
\va{\Control}_{s} \not\in \mathcal{\Control}_{s}\np{ \va{\State}_{s} } } 
\eqfinp 
\end{equation}
Here, \( \GG_{t} \bc{ \sequence{\va{\State}_{s}}{s\in \segment{t}{\horizon}},
\sequence{\va{\Control}_{s}}{s\in \segment{t}{\horizon}} } 
\leq \alpha=1-\beta \) means that 
\( \PP\bc{ \va{\State}_{s} \in \sdo \eqsepv 
\va{\Control}_{s} \in \mathcal{\Control}_{s}\np{ \va{\State}_{s} }
\eqsepv \forall s \geq t } \geq \beta \), for all 
$\PP\in\mathcal{P}$.

\subsubsection{Random exit time and viability}

Let $\mu$ be a measure on the time set~$\TT$ like, for instance, 
the counting measure when $\TT=\NN$ or the Lebesgue 
measure when $\TT=\RR_+$. Then, the random quantity
\begin{equation}
 \mu \{ s \geq t \eqsepv \va{\State}_{s} \not\in \sdo \mtext{ or }
\va{\Control}_{s} \not\in \mathcal{\Control}_{s}\np{ \va{\State}_{s} } \}
\end{equation}
measures the number of times that the state-control path
\( \bp{ \sequence{\va{\State}_{s}}{s\in \segment{t}{\horizon}},
\sequence{\va{\Control}_{s}}{s\in \segment{t}{\horizon}} } \)
\emph{exits} from the viability constraints.

Using risk measures --- like the 
\emph{tail/average/conditional value-at-risk} 
\cite{Follmer-Schied:2002} --- 
or stochastic orders \cite{Muller-Stoyan:2002,Shaked-Shanthikumar:2007},
we have differents ways to express that this random quantity remains ``small''.

\subsubsection{The general umbrella of cost functions}
\label{cost_functions}

All the examples above, and many more \cite{Rouge-Mathias-Deffuant:2015},
fall under the general umbrella of cost functions as follows.

Consider a starting time~$t\in\TT$ 
and a measurable function
\begin{equation}
  \Psi_{t} : \prod_{s=t}^{\horizon}\STATE_{s} \times
  \prod_{s=t}^{\horizon}\CONTROL_{s} \times \Omega \to \RR \eqsepv 
\label{eq:utility}
\end{equation}
that attachs a \emph{disutility} or \emph{cost} ---
the opposite of \emph{value}, \emph{utility}, \emph{payoff} ---
to any tail state and control path, starting from time~$t$,
and to any state of Nature.

Let $\FF$ be a risk measure that maps random variables on~$\Omega$
towards the real numbers. Then, an extended risk measure is given by
\begin{equation}
\GG_{t} \bc{ \sequence{\va{\State}_{s}}{s\in \segment{t}{\horizon}},
\sequence{\va{\Control}_{s}}{s\in \segment{t}{\horizon}} } 
= \FF \bc{ \Psi_{t} 
\bp{ \sequence{\va{\State}_{s}\np{\cdot}}{s\in \segment{t}{\horizon}} ,
\sequence{\va{\Control}_{s}\np{\cdot}}{s\in \segment{t}{\horizon}},
\cdot } } 
\eqfinp
\end{equation}

\subsection{Resilience and risk minimization}

When the set~\( \ResilientStrategies_{t}\np{\bar\state_{t}} \)
of {resilient strategies} at time~$t$ 
in~\S\ref{Resilient_strategies_and_resilient_states} is not empty,
how can we select one among the many? %Or at least restrict the search?
Here is a possible way that makes use of risk measures for 
risk minimization purposes. 

\subsubsection{Indicators of resilience}

Let $\GG_{t}$ be an extended risk measure.
We can look for resilient strategies that minimize risk,
solution of 
\begin{equation}
  \min_{ \strategy \in \ResilientStrategies_{t}\np{\bar\state_{t}} } 
\GG_{t} \bc{
\bp{ \phi_{\segment{t}{\horizon}}^{\strategy}  \np{\bar\state_{t}, \cdot} ,
\cdot } }
\eqfinp 
\end{equation} 
The minimum of the risk measure is a potential candidate
as an \emph{indicator of resilience}.
Indeed, it is the best achievable measure of residual risk 
under a resilient strategy.

\subsubsection{Examples}

Using cost functions as in~\S\ref{cost_functions},
we can look for resilient strategies 
that minimize \emph{expected costs}
\begin{equation}
  \min_{ \strategy \in \ResilientStrategies_{t}\np{\bar\state_{t}} } \EE \bc{
\Psi_{t} \bp{ \phi_{\segment{t}{\horizon}}^{\strategy}  \np{\bar\state_{t}, \cdot} ,
\cdot } } \eqsepv
\end{equation}
or that minimize \emph{worst case costs}, where 
\( \overline\Omega \subset \Omega \), 
\begin{equation}
  \min_{ \strategy \in \ResilientStrategies_{t}\np{\bar\state_{t}} } 
\sup_{ \omega \in \overline\Omega }  
\Psi_{t} \bp{ \phi_{\segment{t}{\horizon}}^{\strategy}  \np{\bar\state_{t}, \omega } ,
\omega  } \eqsepv
\end{equation}
or, more generally, that minimize
\begin{equation}
  \min_{ \strategy \in \ResilientStrategies_{t}\np{\bar\state_{t}} } \FF \bc{
\Psi_{t} \bp{ \phi_{\segment{t}{\horizon}}^{\strategy}  \np{\bar\state_{t}, \cdot} ,
\cdot } } \eqsepv
\end{equation} 
where $\FF$ is a risk measure that maps random variables on~$\Omega$
towards the real numbers \cite{Follmer-Schied:2002}.

For instance, in the robust viability setting of~\S\ref{recovery_time},
an {indicator of resilience} could be the minimum (over all resilient strategies)
of the maximal (over all states of Nature in~$\overline\Omega$) 
recovery time.

\section{Conclusion} 
\label{sec:conclusions}

Resilience is a rehashed concept in natural hazard management.
Most of the formalizations of the concept require that, after any perturbation,
the state of a system returns to an acceptable subset of the state set.
Equipped with tools from control theory under uncertainty, we have proposed
that resilience is the ability for the state-control random process as a whole 
to be driven to an acceptable ``recovery regime'' 
by a proper resilient strategy (adaptive).
Our definition of resilience is a form of controlability:
a resilient strategy has the property to shape the closed loop flow
so that the resulting state and control random process 
belongs to a given subset of random processes, the acceptable recovery regimes.

We have proposed to handle risk thanks to risk measures\footnote{%
We also hinted at the possibility to use so-called stochastic orders.},
by defining recovery regimes that represent a form of ``risk containment''.
In addition, risk measures are potential candidates as 
{indicators of resilience} as they measure the residual risk 
under a resilient strategy.

Our contribution is formal, with its pros and cons:
by its generality, our approach covers a large scope of notions of 
resilience; however, such generality makes it difficult to propose
resolution methods. For instance, the possibility to use dynamic programing
in stochastic viability relies upon a white noise assumption 
that we have not supposed. Much would remain to be done 
regarding applications and numerical implementation.
\bigskip

\paragraph{Acknowledgement.} The author is indebted to the 
editor-in-chief, the advisory editor and two reviewers. 
They supplied detailed critique, comments and inputs which, ultimately, 
contributed to an improved version of the manuscript.

\newcommand{\noopsort}[1]{} \ifx\undefined\allcaps\def\allcaps#1{#1}\fi

\end{document}